\documentclass[12pt]{amsart}
\usepackage{latexsym}
\usepackage{amsmath}
\usepackage{amsfonts}
\usepackage{amsthm}
\usepackage{amssymb}
\usepackage{enumerate}
\usepackage[T1]{fontenc}
\title{On a problem of T. Szostok concerning the Hermite-Hadamard inequalities.}
\newtheorem{theorem}{Theorem}%[section]
\newtheorem{corollary}{Corollary}%[section]
\newtheorem{proposition}{Proposition}%[section]
\theoremstyle{remark}
\theoremstyle{definition}
\newtheorem{definition}{Definition}%[section]
\newcommand{\R}{\mathbb{R}}
\newcommand{\Q}{\mathbb{Q}}

\newcommand{\N}{\mathbb{N}}
\newcommand{\K}{\mathbb{K}}

\begin{document}
\begin{abstract}
In the present paper we solve a problem posed by Tomasz  Szostok who asked about the solutions $f$ and $F$  to the system of inequalities 
$$
f\Big(\frac{x+y}{2}\Big)\leq \frac{F(y)-F(x)}{y-x}\leq \frac{f(x)+f(y)}{2}.
$$
We show that $f$ and $F$ are the solutions to the above system of inequalities if and only if $f$ is a continuous convex function and $F$ is primitive function of $f$. 
This result can be interpreted as a regularity phenomenon-the solutions to the system of functional inequalities turn out to be regular 
without any additional assumptions.
\end{abstract}

\title{On a problem of T. Szostok concerning the Hermite-Hadamard inequalities.}

\author[A. Olbry\'s]{Andrzej Olbry\'s}
\address{Institute of Mathematics, University of Silesia,
40-007 Katowice, ul. Bankowa 14, Poland}
\email{andrzej.olbrys@us.edu.pl}


\maketitle
%\author[A. Olbry\'s›]{Andrzej Olbry\'s›}
%\address{Institute of Mathematics, University of Silesia,
%40-007 Katowice, ul. Bankowa 14, Poland}
%\email{andrzej.olbrys@math.us.edu.pl}

%\subjclass[2000]{Primary yyy, Secondary xxx}
%\keywords{}
%\thanks{This research has been supported by the Hungarian
%Scientific Research Fund (OTKA) Grant K111651.
\section{Introduction}

There are many inequalities valid for convex functions. Probably two of the most well-known ones are the 
Hermite-Hadamard inequalities. They provide a lower and upper estimation for the integral average 
of any convex function defined on a closed interval. More precisely, if $f:[x,y]\to \R$ is a convex function, then 
it is integrable in the sense of Riemann and 
\begin{equation}
f\Big(\frac{x+y}{2}\Big)\leq \frac{1}{y-x}\int_{x}^{y}f(t)dt \leq \frac{f(x)+f(y)}{2}.
\end{equation}

The Hermite-Hadamard inequalities play an important role in research on inequalities and convex analysis. In the literature 
one can find their various generalizations and applications. For more informations on this type of inequalities see the book [7] 
and the references therein. We just note here that first Hermite [13] published these inequalities with some important applications 
and then, ten years later, Hadamard [10] rediscovered their left-hand side. It turns out that in the class of continuous functions 
each of the above inequalities is equivalent to convexity. More general results are given by Rado [23].
Let us observe that we can rewrite the inequalities (1) in the following way 
\begin{equation}
f\Big(\frac{x+y}{2}\Big)\leq \frac{F(y)-F(x)}{y-x}\leq \frac{f(x)+f(y)}{2}, 
\end{equation}
where $F$ is a primitive function of $f$. \\

During the Fourteenth Debrecen-Katowice Winter Seminar on Functional Equations and Inequalities (Hajd\'uszoboszl\'o, January 29 - February 1, 2014) 
Tomasz  Szostok (see [32], Problem 3) asked for  general solutions $f$ and $F$  to the inequalities (2). In the present paper we show that the functions 
$f$ and $F$ satisfy the inequalities (2) if and only if $f$ is continuous and hence convex and $F'=f$. \\ 
This result can be understood as a regularity phenomenon. The solutions to the system of functional inequalities turn out to be regular without any regularity assumptions. In the literature there are known several results of these type for functional
equations but  this phenomenon is very rare for the solutions to the functional
inequalities. As it is well known  the inequality defining convexity, for  functions  defined on an open and convex subset of a   
finite-dimensional real linear space has only continuous solutions ([20], [29]). The same effect was obtained by M. R\v adulescu in the paper [24]. 
 R\v adulescu  showed that if an operator $T:C(X)\to C(X)$ (where $C(X)$ denotes the space of all continuous real valued functions
defined on  compact Hausdorff topological space $X$, equipped with the supremum norm)  satisfies the system of inequalities 
{\setlength{\arraycolsep}{1pt}
$$
\left\{ \begin{array}{ll} 
& T(f+g)\geq T(f)+T(g) \\ 
& T(f\cdot g)\geq T(f)\cdot T(g)
\end{array} \right.  \leqno(6) 
$$}
\hspace{-0,3cm} \noindent for all $f, g\in C(X)$, then  $T$ is linear, multiplicative, and continuous. Z. Ercan in the paper [8] has shown that the assumption that $X$ is a compact Hausdorff space may be dropped.\\

A characterization of functions satisfying the so-called Shannon inequality was obtained by  J. Acz\'el and M. Ostrowski in the paper [3] (see also [2], p. 116). They have shown that the function $h:(0,1)\to \R$ satisfies the functional inequality
$$
\sum_{k=1}^{n}p_k h(p_k)\leq \sum_{k=1}^{n}p_k h(q_k), 
$$
for all such  $p_1,\ldots,p_n;\ q_1,\ldots,q_n \in (0,1)$ that $\sum_{k=1}^{n}p_k= \sum_{k=1}^{n}q_k= 1$ if and only if there exist constants  
$b, c\in \R,\ c\leq 0$ such that 
$$
h(p)=c\log p+b \quad \text{for}\ p\in (0,1).
$$

On the other hand,  W. Fechner [9] considered the two following inequalities similar to ours
$$
f\Big(\frac{x+y}{2}\Big)\leq \frac{f(y)-f(x)}{y-x}\ \ \text{and}\ \ \frac{f(y)-f(x)}{y-x}\leq \frac{f(x)+f(y)}{2},
$$
but in paper [9] the Riemann integrability of solutions was assumed. These two functional inequalities have been also investigated by C. Alsina and R. Ger 
in  [4] who, among other results, gave a general solutions to the first  inequality in the class of  non-negative functions. \\

For functional equations the regularity  phenomenon is characteristic for equations derived from the  mean values theorems.  The first two papers that must be mentioned here are due to J. Acz\'el [1] and S. Haruki [12] who worked with 
$$
F(y)-F(x)=(y-x)f\Big(\frac{x+y}{2}\Big),
$$
and, 
$$
F(y)-F(x)=(y-x)[f(x)+f(y)],
$$
respectively. In this paper we use both expressions considered by Acz\'el and Haruki but with inequalities instead of equalities. 
Results of Acz\'el and Haruki  were generalized among others by M. S. Jacobson, P. Kannappan 
and P. S. Sahoo [14] who dealt with equation 
$$
F(y)-F(x)=(y-x)f(sx+ty).
$$
In all the above mentioned results it turns out that the solution to the considered functional equation  has to be continuous. 
Results in this spirit can be found in the papers [15], [16], [26], [31], [26]  and in monograph [25].
All from the above  equations are a particular case of the general equation of the form  
$$
F(y)-F(x)=(y-x)[f_1(\alpha_1 x+\beta_1 y)+\cdots+f_n(\alpha_n x+\beta_n y)], 
$$
which was considered by T. Szostok in [33]. This functional equation stemming from numerical analysis and is strongly connected with 
classical quadrature rules used in numerical integration.\\

In this paper we obtain the solutions to the system of functional inequalities (2) without any regularity assumptions.  
\vspace{0,5cm}

\section{Main result}
A main tool which we are going to use in the proof of our main result is a concept of so-called the $\K$-Riemann integral. The 
notion of $\K$-Riemann integral was introduced and examined in [21]. Now, we shall use the notation and terminology from [21]. 
We only recall necessary definitions and facts. \\

Let  $\K$ stands for a subfield of the field of real numbers $\R$. Clearly, $\Q\subseteq \K$, where $\Q$ denotes the field of rational numbers. 
 In the sequel the symbol $[a,b]_{A}$ will denote a $A$-convex hull of the set $\{a,b\},$ 
where $A\subseteq \R$ i.e.
$$
[a,b]_{A}=\{\alpha a+(1-\alpha)b: \alpha \in A\cap [0,1]\}.
$$ 
Let $\mathcal{P}_{[a,b]}$ denote the set of partitions of the interval $[a,b]$ i.e. 
$$
\mathcal{P}_{[a,b]}:=\bigcup_{n=1}^{\infty} \{(t_{0},t_{1},...,t_{n}): a=t_{0}<t_{1}<...<t_{n}=b \}.
$$
Following Zs. P\'ales [22] 
we define the set of $\K$-partitions of the interval $[a,b]$ in the following way
$$
{\setlength{\arraycolsep}{1pt}
\begin{array}{ll} 
\mathcal{P}^{\K}_{[a,b]}:&=\Big\{(t_{0},t_{1},...,t_{n})\in \mathcal{P}_{[a,b]}: \frac{t_{i}-a}{b-a}\in \K,\ i=1,2,...,n \Big\} \\
&=\Big\{(t_{0},t_{1},...,t_{n})\in \mathcal{P}_{[a,b]}: t_{i}=a+\alpha_{i}(b-a): \alpha_{i}\in \K\cap [0,1], i=1,2,...,n \Big\}\\
&=\Big\{(t_{0},t_{1},...,t_{n})\in \mathcal{P}_{[a,b]}: t_{i}\in [a,b]_{\K},\ i=1,2,...,n \Big\} 
\end{array}}
$$

Now, suppose that $f:[a,b]\to \R$ is a bounded function on the set $[a,b]_{\K}$ with 
$$
M:=\sup_{x\in [a,b]_{\K}}f(x),\qquad m:=\inf_{x\in [a,b]_{\K}}f(x).
$$
For a given $\K$-partition $\pi=(t_{0},t_{1},...,t_{n}) \in \mathcal{P}^{\K}_{[a,b]}$ let 
$$
M_{i}:=\sup_{x\in [t_{i-1},t_{i}]_{\K}}f(x),\qquad m_{i}:=\inf_{x\in [t_{i-1},t_{i}]_{\K}}f(x),\ \ \ i=1,2,...n.
$$
These suprema and infima are well-defined, finite real numbers since $f$ is bounded on $[a,b]_{\K}$ . Moreover, 
$$
m\leq  m_{i}\leq M_{i}\leq M,\quad i=1,2,...,n.
$$

We define the upper $\K$-Riemann sum of $f$ with respect to the partition $\pi$ by 
$$
U_{\K}(f,\pi):=\sum_{i=1}^{n}M_{i}(t_{i}-t_{i-1}),
$$
and the lower $\K$-Riemann sum of $f$ with respect to the  partition $\pi$ by 
$$
L_{\K}(f,\pi):=\sum_{i=1}^{n}m_{i}(t_{i}-t_{i-1}).
$$
Note that 
$$
m(b-a)\leq L_{\K}(f,\pi)\leq U_{\K}(f,\pi)\leq M(b-a).
$$

Now, we define the upper  $\K$-Riemann integral 
of $f$ on $[a,b]$ by
$$
\overline{\int\limits_{a}}^{b}f(t)d_{\K}t:=\inf\{U_{\K}(f,\pi):\ \pi \in \mathcal{P}^{\K}_{[a,b]}\}
$$
and the lower $\K$-Riemann integral by 
$$
\underline{\int}_{a}^{b}f(t)d_{\K}t:=\sup\{L_{\K}(f,\pi):\ \pi \in \mathcal{P}^{\K}_{[a,b]}\}
$$
\begin{definition}
A function  $f:[a,b]\to \R$  bounded on $[a,b]_{\K}$ is said to be a $\K$-\emph{Riemann integrable} on $[a,b]$ if its 
upper  and lower integral are equal. In that case, the $\K$-Riemann integral of $f$ on $[a,b]$ 
we denote by 
$$
\int\limits_{a}^{b}f(t)d_{\K}t.
$$
\end{definition}

Similar as for usual integral we have the following equivalent condition to $\K$-Riemann integrability.
\begin{theorem}
A function $f:[a,b]\to \R$ is $\K$-Riemann integrable on $[a,b]$ if and only if for every 
sequence $\{\pi_{n}\}_{n\in \N} \subset \mathcal{P}^{\K}_{[a,b]}$, $\pi^{(n)}=(t_{0}^{(n)}, t_{1}^{(n)},...,t_{k_{n}}^{(n)})$ such that 
$$
\max_{1\leq j\leq k_{n}}(t_{j}^{(n)}-t_{j-1}^{(n)})\to _{n\to \infty} 0,
$$  
and for any choice $s_{j}^{(n)}\in [t_{j-1}^{(n)},t_{j}^{(n)}]$ of the partition $\pi_{n}$ we have 
$$
\int_{a}^{b}f(t)d_{\K}t=\lim_{n\to \infty}\sum_{j=1}^{k_{n}}f(s_{j}^{(n)})(t_{j}^{(n)}-t_{j-1}^{(n)}).
$$
\end{theorem}

We have the following evident result: 
\begin{corollary}
 If a function $f:[a,b]\to \R$ is a Riemann integrable in the usual sense, then for an arbitrary field  
 $\K \subseteq \R$ $f$ is a $\K$-Riemann integrable, moreover,
 $$
 \int_{a}^{b}f(t)d_{\K}t=\int_{a}^{b}f(t)dt.
 $$
\end{corollary}

In this place we recall the following well-known definitions: 

\begin{definition}
A mapping $f:\R \to \R$ is called \emph{additive} if it satisfies a Cauchy's functional equation 
$$
f(x+y)=f(x)+f(y),
$$
for every $x, y\in \R$. A map $f$ is called $\K$-\emph{linear} if $f$ is additive and $\K$-homogeneous i.e. 
$$
f(\alpha x)=\alpha f(x),
$$
is fulfilled for every $x\in \R$ and $\alpha \in \K$. 
\end{definition}
It is well-known that every additive function is  $\Q$-homogeneous.
\begin{definition}
A function $f:(a,b)\to \R$ is said to be \emph{Jensen-convex} if 
$$
f\Big(\frac{x+y}{2}\Big)\leq \frac{f(x)+f(y)}{2},
$$
for every $x, y\in (a,b)$. The solution to the corresponding functional equation is said be 
\emph{Jensen-affine function}. It is well-known (see for instance [17]) that any Jensen-affine function $f:(a,b)\to \R$ 
is of the form $$f(x)=a(x)+c,\ x\in (a,b),$$ where $a:\R \to \R$ is an additive function and  $c\in \R$ is a constant.\\

A map $f:(a,b)\to \R$ is called $\K$-convex if 
$$
f(\alpha x+(1-\alpha)y)\leq \alpha f(x)+(1-\alpha)f(y),
$$
for every $x, y\in (a,b)$ and $\alpha \in \K\cap (0,1)$.
\end{definition}
It is known that any Jensen-convex function is also $\Q$-convex (see [13]). On the other hand, if $f$ is $\K$-convex 
then it is also Jensen-convex [4].
\begin{theorem} Any $\K$-convex function $f:[a,b]\to \R$ is a $\K$-Riemann integrable. In particular,  such is any  $\K$-linear function 
$g:\R\to \R$, moreover, 
$$
\int_{a}^{b}g(t)d_{\K}t=g \Big(\frac{a+b}{2}\Big)(b-a).
$$
\end{theorem}

In the proof of our main result we use the following property of $\K$-Riemann integral: 

\begin{proposition}
Let $f, g:[a,b]\to \R$ be  $\K$-Riemann integrable functions. If 
$$
g(x)\leq f(x),\quad x\in [a,b]_{\K},
$$
then 
$$
\int_{a}^{b}g(x)d_{\K}x\leq \int_{a}^{b}f(x)d_{\K}x.
$$ 
\begin{proof}
By the assumption, for arbitrary $\pi \in \mathcal{P}_{[a,b]}^{\K}$ we have 
$$
U_{\K}(\pi,g)\leq U_{\K}(\pi,f),
$$
so 
$$
\int_{a}^{b}g(x)d_{\K}x\leq U_{\K}(\pi,g)\leq U_{\K}(\pi,f),\quad \pi \in \mathcal{P}_{[a,b]}^{\K}.
$$
Taking an infimum over the all $\pi \in \mathcal{P}_{[a,b]}^{\K}$ we get 
$$
\int_{a}^{b}g(x)d_{\K}x\leq \overline{\int_{a}^{b}}f(x)d_{\K}x= \int_{a}^{b}f(x)d_{\K}x,
$$
which finishes the proof.
\end{proof}
\end{proposition}

Now, we are able to proof our main result. The following theorem gives the general solutions to the system of inequalities (2).
\begin{theorem} Let $f, F:(a,b)\to \R$ be  given functions. Then they satisfy the inequalities 
$$
f\Big(\frac{x+y}{2}\Big)\leq \frac{F(y)-F(x)}{y-x}\leq \frac{f(x)+f(y)}{2},\quad x, y\in (a,b),\ x\neq y, 
$$
if and only if $f$ is  convex function and $F$ is  primitive function of $f$.
\begin{proof} If $f$ is  convex function and $F$ is primitive of $f$ then the inequalities (2) hold on account of 
the classical Hermite-Hadamard inequalities. Conversely, assume that $f$ and $F$ satisfy the inequality (2), in particular,  
$f$ is convex in the sense of Jensen. Fix $x, y\in (a,b),\ x\neq y$, say $x<y$. 
Let $$\pi_n:=(x_{0}^{(n)}, x_{1}^{(n)},\dots, x_{k_n}^{(n)})\in \mathcal{P}^{\Q}_{[x,y]},\quad n\in \N,$$ be an arbitrary 
sequence of $\Q$-partitions of the interval $[x,y]$ such that
$$
\max_{1\leq j\leq k_{n}}(x_{j}^{(n)}-x_{j-1}^{(n)})\to _{n\to \infty} 0.
$$  
By (2) 
$$
f\Big(\frac{x_{j}^{(n)} +x_{j-1}^{(n)}}{2}\Big)\leq \frac{F(x_{j}^{(n)})-F(x_{j-1}^{(n)})}{x_{j}^{(n)}-x_{j-1}^{(n)}}\leq \frac{f(x_{j}^{(n)}) 
+f(x_{j-1}^{(n)})}{2},\quad j=1,\dots, k_n,
$$
or equivalently, for $j=1,\dots, k_n $ we have
$$
f\Big(\frac{x_{j}^{(n)} +x_{j-1}^{(n)}}{2}\Big)(x_{j}^{(n)}-x_{j-1}^{(n)})\leq F(x_{j}^{(n)})-F(x_{j-1}^{(n)})\leq \frac{f(x_{j}^{(n)}) 
+f(x_{j-1}^{(n)}}{2}(x_{j}^{(n)}-x_{j-1}^{(n)}).
$$
Summing up the above inequalities from $j=1$ to $k_n$ we obtain

{\setlength{\arraycolsep}{1pt}
\begin{eqnarray*}
\sum_{j=1}^{k_n}f \Big(\frac{x_j^{(n)}+x_{j-1}^{(n)}}{2} \Big) &(& x_{j}^{(n)}-x_{j-1}^{(n)} ) \leq  \sum_{j=1}^{k_n}(F(x_{j}^{(n)})-F(x_{j-1}^{(n)})  \nonumber \\
 &\leq &  \frac{1}{2}\Big[\sum_{j=1}^{k_n}f(x_{j}^{(n)})(x_{j}^{(n)}-x_{j-1}^{(n)}) + \sum_{j=1}^{k_n}f(x_{j-1}^{(n)})(x_{j}^{(n)}-x_{j-1}^{(n)})\Big] \nonumber
\end{eqnarray*}}

On account of Theorem 6 $f$ is a $\Q$-Riemann integrable, so taking the limit as $n\to \infty$, because $$\sum_{j=1}^{k_n}(F(x_{j}^{(n)})-F(x_{j-1}^{(n)})) =F(y)-F(x),\quad n\in \N,$$  we get
$$
\int_{x}^{y}f(t)d_{\Q}t\leq F(y)-F(x)\leq \frac{1}{2}\Big(\int_{x}^{y}f(t)d_{\Q}t+\int_{x}^{y}f(t)d_{\Q}t\Big),
$$ 
therefore, 
$$
F(y)-F(x)=\int_{x}^{y}f(t)d_{\Q}t.
$$

Note that, this property implies that the $\Q$-Riemann integral of $f$ is additive with respect to the interval of integration. 
(in general it is not true see [21, Remark 19]). Indeed, for arbitrary $\alpha< \gamma < \beta,\ \alpha, \beta \in (a,b)$ we have 
{\setlength{\arraycolsep}{1pt}
\begin{eqnarray*}
\int_{\alpha}^{\beta}f(t)d_{\Q}t&=&F(\beta)-F(\alpha)=F(\beta)-F(\gamma)+F(\gamma)-F(\alpha) \nonumber \\ 
&=&  \int_{\gamma}^{\beta}f(t)d_{\Q}t +  \int_{\alpha}^{\gamma}f(t)d_{\Q}t.  \nonumber
\end{eqnarray*}}

Now, we show that $f$ is a convex function. To do it, fix $x, y\in (a,b),\ \lambda \in (0,1)$ arbitrary and put 
$z:=\lambda x+(1-\lambda)y$. By Rod\'e's theorem (see [30]) there exists a Jensen-affine support function at the point $z$  i.e. a Jensen-affine function 
$g_z:(a,b)\to \R$ such that $g_z(z)=f(z)$,  and 
 $$g_z(x)\leq f(x),\ x\in (a,b).$$
Since $g_z$ has the form 
$$
g_z(x)=a(x)+c\quad \text{for}\ x\in (a,b), 
$$
where $a:\R\to \R$ is an additive function and $c\in \R$ is a constant, then it follows from Corollary 3   and Theorem 6 that for $\alpha, \beta\in (a,b),\ 
\alpha<\beta$ we have 
$$
\int_{\alpha}^{\beta}g_z(x)d_{\Q}x=a\Big(\frac{\alpha+\beta}{2}\Big)(\beta-\alpha)+c(\beta-\alpha).
$$
Using this and  on account of Proposition 7 for arbitrary $\gamma \in (\alpha,\beta)$ we get 
{\setlength{\arraycolsep}{1pt}
\begin{eqnarray}
\frac{1}{2}a(\gamma)(\beta-\alpha)+ \frac{1}{2}a(\alpha-\beta)\gamma &-& \frac{1}{2}a(\alpha)\alpha+\frac{1}{2}a(\beta)\beta+c(\beta-\alpha)   \nonumber \\ 
&=& a\Big(\frac{\alpha+\gamma}{2}\Big)(\gamma-\alpha) +c(\gamma - \alpha)+a\Big(\frac{\gamma+\beta}{2}\Big)(\beta-\gamma) +c(\beta-\gamma) \nonumber \\ 
&=& \int_{\alpha}^{\gamma}g_{z}(x)d_{\Q}x+\int_{\gamma}^{\beta}g_{z}(x)d_{\Q}x \nonumber \\ 
&\leq & \int_{\alpha}^{\gamma}f(x)d_{\Q}x+\int_{\gamma}^{\beta}f(x)d_{\Q}x=\int_{\alpha}^{\beta}f(x)d_{\Q}x<\infty. \nonumber
\end{eqnarray}}

Let us define a function $h:(a,b)\to \R$ by the formula 
$$
h(\gamma):=\frac{1}{2}a(\gamma)(\beta-\alpha)+\frac{1}{2}a(\alpha-\beta)\gamma-\frac{1}{2}a(\alpha)\alpha+\frac{1}{2}a(\beta)\beta+c(\beta-\alpha). 
$$ 
It is easy to see that $h$ is  Jensen-affine  and, as we have already seen, bounded from above on a non-empty and open subset of $\R$, 
whence on account of the famous Bernstein and Doetsch  theorem [5] continuous. In particular, $a$ is a continuous additive function, therefore 
$$
g_{z}(tu+(1-t)v)=tg_{z}(u)+(1-t)g_{z}(v)\quad \text{for} \ u, v \in (a,b),\ t\in (0,1).
$$
From this we obtain 
$$
f(\lambda x+(1-\lambda)y)=f(z)=g_z(z)=\lambda g_z(x)+(1-\lambda)g_z(y)\leq \lambda f(x)+(1-\lambda)f(y).
$$
Due to the arbitrariness of $x, y\in (a,b)$ and $\lambda\in [0,1]$ we infer that $f$ is  convex function.\\

Now, we show that $F$ is primitive function of $f$. Indeed, for arbitrary $x\in(a,b)$  by putting $y:=x+h$ in (2), (where 
$h\in\R$ is such that $x+h\in (a,b)$) we obtain
$$
f\Big(x+\frac{h}{2}\Big)\leq \frac{F(x+h)-F(x)}{h}\leq \frac{f(x+h)+f(x)}{2}. 
$$
Since any convex function defining on an open and convex subset of a finite dimensional space has to be continuous, then 
letting $h\to 0$ we get
$$
f(x)=F'(x),\quad x\in (a,b).
$$
The proof of our theorem is finished.
\end{proof}

\end{theorem}

\end{document}